\documentclass[12pt]{amsart}
\usepackage{amssymb}
\begin{document}

\title{BRAID GROUP OF A GENETIC CODE}

\author[Bardakov]{Valerij G. Bardakov}
\address{Sobolev Institute of Mathematics, Novosibirsk 630090, Russia}
\email{bardakov@math.nsc.ru}

%
%
\date{\today}

\begin{abstract}
For every group genetic code with finite number of generating and at most with one defining relation
we introduce the braid group of this genetic code. This construction includes the braid group
of Euclidean plane, the braid groups of closed  orientable surfaces, $B$ type groups of
Artin-Brieskorn, and allow us to study all these groups with common point of view. Then we clarify a structure
of the braid group of genetic code, construct a normal form of words, investigate a torsion,
calculate the width of verbal subgroups and establish some other properties. Also we prove that
Scott's system of defining relations for the braid groups of closed surfaces is contradictory.

\end{abstract}

\maketitle


In last years appeared many articles dealing with the braid
groups and  their different
generalizations \cite{BS}-\cite{Bar-F04}. This fact is explained
by connection of theory of braid groups with knot theory. Knot
theory draws attention of many researchers with connection to
appearance of new methods and directions, which arise in this
theory. It is enough to remember about results of Jones \cite{J},
who constructed new knot polynomial invariant, about results of
V.~A.~Vasilieva \cite{Vas}, who constructed  finite type
invariants, about results of Bigelow and Krammer \cite{Big, Kra},
who proved that the braid groups are linear, and also about
result of Dehornoy \cite{D} which said that the braid groups are
right orderable. Apart from that  knot theory is fast connects
with other areas of mathematics and physics \cite{A1, A2}.

One of  the first generalizations of the braid groups are braid groups of closed orientable surfaces,
which were introduced by Zariski \cite{Z}. Then in the Scott's article  \cite{S} were considered
the braid groups of closed orientable surfaces and were found  generators and defining
relations of these groups and also of pure braid groups of these surfaces (a set of
generators and defining relations is called the {\it genetic code} of this group). Unfortunately the
set of defining relations which was written in the Scott's article  is contradictory (see \S~2 of the
present paper). Braid groups of sphere and projecting plane differ from the braid groups of other
surfaces, and as a rule are studied separately (see \cite{GilB, FadB, VB}).

Gonz\'{a}lez-Meneses \cite{Gon} found a new presentation of the braid groups of closed surfaces
in the form
of a set of generators and defining relations. Then Bellingeri \cite{Bel} simplified this presentation and
found generators and defining relations of the braid groups of surfaces with boundary. In preprint
\cite{BG} was found a positive presentation of the braid groups of orientable surfases.
Apart from that,
in \cite{Gon, Bel}
was solved the word problem for corresponding braid groups, and in preprint \cite{BG} was solved
the conjugacy problem for the braid group $B_2(T^2)$ on two strands of torus $T^2$. In general case
the conjugacy problem for the braid groups of surfaces is still open.

In the present paper we recall a definition of the braid groups of manifold (see \S~1), show
that Scott's presentation is contradictory (see \S~2). As Bellingery stated
to the author that
last fact was also observed by V.~S.~Kulikov and Shimada. Then in \S~3 for each genetic code with
finite number of generators and at most one relation we  introduce the braid group of
this genetic code. This constriction includes the braid group of plane,
the braid groups of closed orientable surfaces, the $B$-type of braid group of Artin-Briskorn,
and it allows us to study all these groups with unified position. In \S~4 we
establish a structure
of the braid group of genetic code, construct a normal form of words, investigate torsion,
calculate width of verbal subgroups and establish some others properties.

The results of the present paper were announced in 2000 (see \cite{Bar-T00}), but at the series of
reasons wasn't  published in proper time.

\subsection*{Acknowledgments}
I am very grateful to the participants of the
seminar ``Evariste Galois'' who listened proofs and forwarded series of useful remarks and
suggestions. Separately I thank  M.~V.~Neshchadim for fruitful discussions and useful
conversations. Also I thank
 V.~V.~Vershinin, who introduced  Zariski's article to me. Special thanks goes to my brother
Victor Bardakov who helped me to translate this article into English.

\vskip 20pt

\centerline{ \bf \S~1. Braid groups of manifolds}

\vskip 8pt

We recall  some known facts about the braid groups of manifolds (see \cite{Bir1, Lin}).

Let $M$ be a connected manifold of dimension $\geq 2$. {\it Configuration space} $F_n M$,
$n\in \mathbb{N}$, of manifold $M$ is called the manifold
$$
F_n M=\{ (z_1, z_2, \ldots, z_n)\in \prod^n_{i=1} M~ \vert ~z_i\neq z_j~~
\mbox{for}~~i\neq j \}
$$
of ordering collections of $n$ distinct points from $M$. Its fundamental group
$\pi_1(F_n M)$ is called the $n${\it --strand pure braid group of manifold} $M$ and is denoted
by $P_n(M)$. On the set $F_n M$ is defined an equivalence relations, $z \sim z'$ if the vector
$z = (z_1, z_2, \ldots, z_n) \in F_n M$ different from the vector $z' = (z_1', z_2', \ldots, z_n')
\in F_n M$ by permutation components.
 The fundamental group $\pi_1(F_n M / _{\sim})$ of the factor manifold $F_n M / _{\sim}$ is called the
{\it $n$-strand braid group of manifold} $M$ and denoted by $B_n(M)$.
The natural projection $p : F_n M \longrightarrow F_n M /_{\sim}$ is a regular covering projection.
A group of covering transformations is the symmetric group $S_n$. Hence, there is a canonical
isomorphism
$$
\pi_1(F_n M /_{\sim}) / \pi_1(F_n M) \simeq B_n(M) / P_n(M) \simeq S_n.
$$
If $M$ is a closed, smooth manifold then the map $F_n M \longrightarrow \prod_{i=1}^n M$
induces the epimorphism $P_n(M)\longrightarrow \pi_1(M)\times \ldots \times \pi_1(M)$ of
the pure braid group $P_n(M)$ on
 direct product $n$ copies of fundamental group $\pi_1(M)$.
And what is more, if the dimension of $M$ is greater than two then this epimorphism is injective.
Hence, the braid groups
of manifolds of dimension 2 represent
the largest interest. Further we will assume that $M$ is a two dimension manifold.
Braid group $B_n=B_n(E^2)$ of Euclidean
plane $E^2$ is
defined  by generators
$\sigma_1,\sigma_2,...,\sigma_{n-1}$ and defining relations
$$
\sigma_i\sigma_{i+1}\sigma_i=\sigma_{i+1}\sigma_i\sigma_{i+1},~~\mbox{if}~~~ i=1,2,...,n-2,
$$
$$
\sigma_i \sigma_j = \sigma_j \sigma_i,~~\mbox{if}~~ ~|i-j|\geq 2.
$$

The group $P_n=P_n(E^2)$ is generated by elements $a_{ij},$ $1 \leq i < j \leq n$,
which express through the generators of $B_n$ by the following manner:
$$
a_{i,i+1}=\sigma_i^2,
$$
$$
a_{ij}=\sigma_{j-1}\sigma_{j-2}\ldots\sigma_{i+1}\sigma_i^2\sigma_{i+1}^{-1}\ldots
\sigma_{j-2}^{-1}\sigma_{j-1}^{-1},~~~i+1< j \leq n.
$$
For the pure braid group $P_n$ there is an epimorphism $\eta : P_n \longrightarrow P_{n-1}$,
its kernel $U_n = \mbox{ker}(\eta)$ is a free group of rank $n-1$.
The pure braid group $P_n$ is a
semi direct product of normal subgroup $U_n$, which is free generated by
$a_{1,n},$ $a_{2,n},\ldots,a_{n-1,n}$ and the group $P_{n-1}$.
Similar, $P_{n-1}$ is a semi direct product of the free group $U_{n-1}$ with free generators
$a_{1,n-1},$ $a_{2,n-1},\ldots,a_{n-2,n-1}$ and the subgroup $P_{n-2}$, and so on.
Hence, the group $P_n$
has the following decomposition
$$
P_n=U_n\leftthreetimes (U_{n-1}\leftthreetimes (\ldots \leftthreetimes
(U_3\leftthreetimes U_2))\ldots),
~~~ U_i \simeq F_{i-1}, ~~~i=2, 3, \ldots, n,
$$

The pure braid group $P_n$ is defined by
the relations (for $\varepsilon = \pm$)
$$
\begin{array}{lll}
1) & a_{ik}^{-\varepsilon }\, a_{kj} \, a_{ik}^{\varepsilon } \, = \, \left( a_{ij} \,
a_{kj} \, \right)^{\varepsilon } \, a_{kj} \,
\left( a_{ij} \, a_{kj}\right) ^{-\varepsilon}, & \\
& \\
2) & a_{km}^{-\varepsilon } \, a_{kj} \, a_{km}^{\varepsilon } \, = \, \left( a_{kj} \,
a_{mj}\right) ^{\varepsilon } \, a_{kj} \,
\left( a_{kj} \, a_{mj}\right) ^{-\varepsilon }, &~~~m < j,\\
& \\
3) & a_{im}^{-\varepsilon } \, a_{kj} \, a_{im}^{\varepsilon } \, = \,
\left[ a_{ij}^{-\varepsilon }, \, a_{mj}^{-\varepsilon }\right] ^{\varepsilon } \, a_{kj} \,
\left[ a_{ij}^{-\varepsilon }, \, a_{mj}^{-\varepsilon }\right] ^{-\varepsilon }, & ~~~i < k < m,\\
& \\
4) & a_{im}^{-\varepsilon} \, a_{kj} \, a_{im}^{\varepsilon } \, = \, a_{kj}, & ~~~k < i;~~~m < j~~
\mbox{or}~~ m < k,
\end{array}
$$
where $[a, b] = a^{-1} b^{-1} a b$.

In this case generators of the braid group act on generators of the pure braid group by
the formulas

$$
\begin{array}{lll}
1) & \sigma_{k}^{-\varepsilon }a_{ij}\sigma_{k}^{\varepsilon }=a_{ij}, &
k\neq i-1, i, j-1, j,\\
& \\
2) & \sigma_{i}^{-\varepsilon }a_{i,i+1}\sigma_{i}^{\varepsilon }=a_{i,i+1}, & \\
& \\
3) & \sigma_{i-1}^{-1}a_{ij}\sigma_{i-1}=a_{i-1,j}, & \\
& \\
4) & \sigma_{i-1} a_{ij} \sigma_{i-1}^{-1}=a_{i-1,i} a_{i-1,j}
a_{i-1,i}^{-1}, & \\
& \\
5) & \sigma_{i}^{-1} a_{ij} \sigma_{i}=a_{i+1,j} [ a_{i,i+1}^{-1},
a_{ij}^{-1} ], & j\neq i+1, \\
& \\
6) & \sigma_{i} a_{ij} \sigma_{i}^{-1}=a_{i+1,j}, & j\neq i+1, \\
& \\
7) & \sigma_{j-1}^{-1}a_{ij}\sigma_{j-1}=a_{i,j-1}, & \\
& \\
8) & \sigma_{j-1} a_{ij} \sigma_{j-1}^{-1}=a_{ij}^{-1} a_{i,j-1}
a_{ij}, & \\
& \\
9) & \sigma_{j}^{-1} a_{ij} \sigma_{j}=a_{ij} a_{i,j+1}
a_{ij}^{-1}, & \\
& \\
10) & \sigma_{j} a_{ij} \sigma_{j}^{-1}=a_{i,j+1}, & 1 \leq i < j \leq n-1,\\
\end{array}
$$
(see \cite{Mar}).

If $n \geq 3$ then the center of the group $B_n$ is an infinite cyclic group,
generated by an element
$$
(\sigma_1 \sigma_2 \ldots \sigma_{n-1})^n = a_{12} (a_{13} a_{23})
\ldots (a_{1n} a_{2n} \ldots a_{n-1,n}).
$$

If $M$ is a compact two dimension manifold which differs from the sphere
$S^2$ and projective plane $P^2$ then the braid group $B_n(E^2)$ of Euclidean plane
$E^2$ includes into the braid group $B_n(M)$ of manifold $M$ for all $n \geq 1$.

Fadell and Noivirt proved that if  $M$ is a compact two dimension manifold which differs from
$S^2$ and $P^2$ then there are not elements of finite order in the braid group  $B_n(M)$.

The braid group $B_n(S^2)$  of sphere $S^2$ is generated by elements
$\delta_1, \delta_2,\ldots,\delta_{n-1}$ and defined by relations
$$
\delta_i \, \delta_{i+1} \, \delta_{i} \, = \, \delta_{i+1} \, \delta_{i} \, \delta_{i+1}~~~\mbox{if}~~~ i=1, 2, \ldots ,n-2,
$$
$$
\delta_i \, \delta_{j} \, = \, \delta_{j} \, \delta_{i}~~~\mbox{if}~~~ |i-j| \geq 2,
$$
$$
\delta_1\, \delta_{2}\ldots \delta_{n-2} \, \delta_{n-1}^2 \,\delta_{n-2}\ldots
\delta_{2} \, \delta_1 \, = \, 1.
$$
In this case $P_2(S^2) = 1$, $P_3(S^2) = B_2(S^2) \simeq \mathbb{Z}_2$, $B_3(S^2)$
is a meta-cyclic group of order 12. If
$n \geq 3$ then the center of $B_n(S^2)$ is a subgroup of order 2,  generated
by an element $(\delta_1 \delta_2 \ldots \delta_{n-1})^n$. The group $B_{n-1}(S^2)$ is not included
in the group $B_n(S^2)$ if $n \geq 3$.

\vskip 20pt

\centerline{\bf \S~2.  On Scott's presentation}


\vskip 5pt

If $T_g$ is a compact orientable surface of genus $g \geq 1$ then generators and relations of
the braid group  $B_n(T_g)$, $n\geq 1$,
was found in the work \cite{S}. We will show that this system of defining relations
is contradictorily.

Let $T=T_1$ be a compact orientable surface of genus 1 (torus).
Let (US) consider the braid group $B_2(T)$ on two strands and its subgroup
of the pure braid
$P_2(T)$. As follows from Theorem 1.3 of Scott's paper, the group $P_2(T)$
is generated by elements
$$
a_{12},~~~\rho_{11},~~~\rho_{12},~~~\rho_{21},~~~\rho_{22}
$$
and defined by relations
$$
\begin{array}{llr}
[\rho_{i2}, \rho_{i1}] = a_{12}, &
i=1,2, & ~~~~~~~~~~(1.1)\\
& \\
\rho_{2i}^{-1} \rho_{1i} \rho_{2i} = a_{12}^{-1} \rho_{1i} a_{12}, &
i=1,2, & ~~~~~~~~~~(2.1)\\
& \\
\rho_{2i} \rho_{1i} \rho_{2i}^{-1} = \rho_{1i}^{-1} a_{12} \rho_{1i} a_{12}^{-1} \rho_{1i}, &
i=1,2, & ~~~~~~~~~~(2.2)\\
& \\
\rho_{1i}^{-1} \rho_{2i} \rho_{1i} = \rho_{2i} a_{12}^{-1} \rho_{2i} a_{12} \rho_{2i}^{-1}, &
i=1,2, & ~~~~~~~~~~(2.3)\\
& \\
\rho_{1i} \rho_{2i} \rho_{1i}^{-1} =  a_{12} \rho_{2i} a_{12}^{-1}, &
i=1,2, & ~~~~~~~~~~(2.4)\\
& \\
\rho_{21}^{-1} \rho_{12} \rho_{21} = a_{12}^{-1} \rho_{11} a_{12} \rho_{11}^{-1} \rho_{12} a_{12}, &
 & ~~~~~~~~~~(3.1)\\
& \\
\rho_{21} \rho_{12} \rho_{21}^{-1} = \rho_{11}^{-1} a_{12} \rho_{11} a_{12}^{-1} \rho_{12}
 \rho_{11}^{-1} a_{12}^{-1} \rho_{11}, &
 & ~~~~~~~~~~(3.2)\\
& \\
\end{array}
$$
$$
\begin{array}{llr}
\rho_{11}^{-1} \rho_{22} \rho_{11} = \rho_{22} \rho_{21} a_{12} \rho_{21}^{-1}, &
 & ~~~~~~~~~~(3.3)\\
& \\
\rho_{11} \rho_{22} \rho_{11}^{-1} =  \rho_{22} a_{12}^{-1}, &
 & ~~~~~~~~~~(3.4)\\
& \\
\rho_{22}^{-1} \rho_{11} \rho_{22} = a_{12}^{-1} \rho_{11}, &
 & ~~~~~~~~~~(4.1)\\
& \\
\rho_{22} \rho_{11} \rho_{22}^{-1} = \rho_{12}^{-1} a_{12} \rho_{12} \rho_{11}, &
 & ~~~~~~~~~~(4.2)\\
& \\
\rho_{12}^{-1} \rho_{21} \rho_{12} = \rho_{22} a_{12}^{-1} \rho_{22}^{-1} \rho_{21} a_{12}^{-1}
\rho_{22} a_{12} \rho_{22}^{-1}, &
 & ~~~~~~~~~~(4.3)\\
& \\
\rho_{12} \rho_{21} \rho_{12}^{-1} = a_{12} \rho_{21} \rho_{22}^{-1} a_{12} \rho_{22} a_{12}^{-1}. &
 & ~~~~~~~~~~(4.4)\\
& \\
\end{array}
$$
From Theorem 1.4 follows that the group $B_2(T)$ is generated by elements
$$
\sigma_1,~~~\rho_{11},~~~\rho_{12},
$$
which are connected with generators of group $P_2(T)$ by equations
$$
a_{12} = \sigma_1^2,~~~\rho_{21} = \sigma_1^{-1} \rho_{11}
\sigma_1,~~~\rho_{22} = \sigma_1^{-1} \rho_{12}
\sigma_1.
$$
A group $A_1 = \langle a_{12}, \rho_{11}, \rho_{12} \rangle $ is isomorphic
to a free group with free generators
$\rho_{11}, \rho_{12}$, which is normal in the group $P_2(T)$.

From the relation (1.1) follows that in the group $A_1$ realizes the relation
$$
[\rho_{12}, \rho_{11}] = a_{12}. \eqno{(5)}
$$
There conjugate both hand parts of this relation by element $\rho_{22}$. Then in the left
hand part using relation (2.1) and (4.1) we get the equation
$$
[\rho_{12}, \rho_{11}]^{\rho_{22}} = [a_{12}^{-1} \rho_{12} a_{12}, a_{12}^{-1}
\rho_{11}] = a_{12}^{-1} \rho_{12}^{-1} a_{12} \rho_{11}^{-1}
\rho_{12} \rho_{11}.
$$
In the right hand part we get the equation
$$
\rho_{22}^{-1} a_{12} \rho_{22} = a_{12}^{-1} \rho_{12} a_{12} \rho_{12}^{-1}
a_{12}. \eqno{(6)}
$$
Really, from the relation (2.4) we have the equation
$$
\rho_{22}^{-1} a_{12} \rho_{22} = (\rho_{22}^{-1} \rho_{12} \rho_{22})
\rho_{12}^{-1} a_{12}.
$$
Applying to the expression in the brackets relation (2.1) unde
$i=2$ we get the equation (6). Therefore, in the free group $A_1$ we
have the equation
$$
a_{12}^{-1} \rho_{12}^{-1} a_{12} \rho_{11}^{-1} \rho_{12} \rho_{11}=
a_{12}^{-1} \rho_{12} a_{12} \rho_{12}^{-1}
a_{12},
$$
or, after multiplication on the left on  $a_{12}$:
$$
\rho_{12}^{-1} a_{12} \rho_{11}^{-1} \rho_{12} \rho_{11}=
\rho_{12} a_{12} \rho_{12}^{-1}
a_{12}. \eqno{(7)}
$$
The group $A_1$ is free generated by elements $\rho_{11}, \rho_{12}$, and element $a_{12}$
is expressed through these elements with help of equation (5).
There escape in the equation
(7) from the generator  $a_{12}$ we get the equation
$$
\rho_{12}^{-2} \rho_{11}^{-1} \rho_{12} = \rho_{11}^{-1}
\rho_{12} \rho_{11} \rho_{12}^{-2}
\rho_{11}^{-1}.
$$
But in the free group with free generators $\rho_{11}$ and $\rho_{12}$ such equation
impossible. Hence, the system of defining relations of Scott for the braid groups of two
dimension orientable manifolds is contradiction. Similar, there can show
that although  the system  of defining relations which is constructed in Theorems 1.1 и 1.2
for the braid groups of non-orientable manifolds is contradiction too.

\vskip 20pt

\begin{center}

{\bf \S~3.  Braid group of genetic code }

\end{center}

\vskip 5pt

Let us consider a genetic code
$$
\mathcal{G} = \langle X \parallel R \rangle, \eqno{(1)}
$$
where $X=\{ x_1,\ldots, x_l \}$ is a finite set of letters and $R$
is an empty set or a set which contains one word in the alphabet
$X^{\pm 1}=\{ x^{\pm 1}_1,\ldots, x^{\pm 1}_l \}$. At that we will
distinguish cases when  $R=\varnothing$ is an empty set and
$R=\{ e \}$ is a set from the trivial word. The group which is defined by the genetic code (1)
we will denote by $G(\mathcal{G})$. We compare to the genetic code  (1) and a natural number
$n$ a group $B_n(\mathcal{G})$ which will call
$n$---{\it strand braid group of genetic code } $\mathcal{G}$. If $n=1$ then we let
$B_1(\mathcal{G})=G(\mathcal{G})$. If  $n > 1$, $R=\{ r \}$ then the group
$B_n(\mathcal{G})$ is defined by generators $\sigma_1,\sigma_2,...,\sigma_{n-1},$
$x_1,\ldots, x_l$ and defining relations

$$
\sigma_i \, \sigma_{i+1} \, \sigma_i \, = \, \sigma_{i+1} \, \sigma_i \, \sigma_{i+1},~~~
i=1,2,...,n-2, \eqno{(2)}
$$

$$
\sigma_i  \, \sigma_j \,  = \,  \sigma_j  \, \sigma_i,~~~ |i-j|\geq 2,
\eqno{(3)}
$$

$$
\sigma_i  \, x_j \,  = \,  x_j \,  \sigma_i,~~~ i=2, 3,..., n-1,~~~j=1, 2,...,l,
\eqno{(4)}
$$

$$
x_j (\sigma_1^{-1} \,  x_j  \, \sigma_1^{-1})  \, =  \, (\sigma_{1}^{-1} \,  x_j \,  \sigma_1^{-1})
 \,  x_{j},
~~~ j=1, 2,..., l, \eqno{(5)}
$$

$$
x_j (\sigma_1^{-1} \,  x_i^{-1} \,  \sigma_1) \,  x_j^{-1} \,   (\sigma_{1}^{-1} \,  x_i  \,
\sigma_1)  \, = \,  \sigma_{1}^2,
~~~ 1\leq i < j \leq l, \eqno{(6)}
$$

$$
r = \sigma_1  \, \sigma_2 \,  \ldots \,  \sigma_{n-2} \,  \sigma_{n-1}^{2}  \, \sigma_{n-2} \,
 \ldots \,  \sigma_2 \,
\sigma_{1}.
 \eqno{(7)}\\
$$

Such defined class groups $B_n(\mathcal{G})$ includes many known groups.
Let us give a few examples.

If $\mathcal{G}_{\varnothing}^{\varnothing} = \langle \varnothing \parallel \varnothing \rangle $
is the empty genetic code then there are not exist the relations
(4)--(7) and $B_n(\mathcal{G}_{\varnothing}^{\varnothing})=B_n(E^2)$ is the braid group
of Euclidean plane. If $\mathcal{G}_{0}^{\varnothing} = \langle \varnothing \parallel e=e \rangle $,
then there are not exist the relations  (4)--(6) and
$B_n(\mathcal{G}_0^{\varnothing})=B_n(S^2)$ is the braid group of two dimension sphere.

Consider the genetic codes
$$
\mathcal{G}^l_{\varnothing} = \langle x_1, x_2, \ldots, x_l ~ \vert \vert ~ \varnothing \rangle,~~~
\mathcal{G}^l_{0} = \langle x_1, x_2, \ldots, x_l ~ \vert \vert ~ e = e \rangle,
$$
$$
\mathcal{G}^l_{1} = \langle x_1, x_2, \ldots, x_l ~ \vert \vert ~ r = e \rangle,~~~l \geq 1,
$$
where $r$ is a non-trivial word.
At that we consider that the group $G(\mathcal{G}^{l}_{1})$ is not a free group
since in the opposite case we can define (IT) by the genetic code with the empty set of relations.
Evidently, $G(\mathcal{G}^l_{\varnothing}) \simeq G(\mathcal{G}^l_{0}) \simeq F_l$
is the free group with free generators $x_1, x_2, \ldots, x_l$;
$G(\mathcal{G}^l_{1})$ is the group with one defining relation, and this group is
a homomorphic image of groups
 $G(\mathcal{G}^l_{\varnothing})$
and $G(\mathcal{G}^l_{0})$.

If
$$
\mathcal{G}_1^{2p}= \langle x_1, \ldots, x_{2p} \parallel x_1x_2^{-1}x_3 \ldots x_{2p}^{-1}x_1^{-1}x_2x_3^{-1}
\ldots x_{2p} \rangle, ~~~p \geq 1,
$$
is the genetic code of fundamental group of compact orientable surface  genus
$p$ then $B_n(\mathcal{G}_1^{2p})$ is the braid group of this surface, it was defined by
Zariski \cite{Z}. If $\mathcal{G} = \langle x \parallel \varnothing \rangle$ then
$B_n(\mathcal{G})$ is an Artin braid group corresponding to the Coxeter group of type $B$
\cite{Bur}. Indeed,
Coxeter group of type $B$ corresponds to the Artin group with generators
 $\sigma_1,\sigma_2,...,\sigma_{n-1},$ $\tau$, where the generators
$\sigma_1,\sigma_2,...,\sigma_{n-1}$ satisfy to relations
(2)--(3), and the relations which contain $\tau $ have the form

$$
\sigma_i \tau = \tau \sigma_{i},~~~ i=2,3,...,n-1, \eqno{(8)}
$$

$$
\tau \sigma_1 \tau \sigma_1 = \sigma_1 \tau \sigma_1 \tau, \eqno{(9)}
$$
i. e. relations (8) correspond to relations (4), and rewriting the relations (9) in the form
$$
\tau^{-1} \sigma_1^{-1} \tau^{-1} \sigma_1^{-1} = \sigma_1^{-1} \tau^{-1} \sigma_1^{-1} \tau^{-1},
$$
and put $\tau^{-1} = x$ we get the braid group $B_n (\langle x \parallel \varnothing \rangle )$.
Note that this group connects to links in body torus similar as the
braid group of Euclidean plane connects to links in the three dimension sphere
$S^3$ (see \cite{Ver}).

The braid group  $B_n(\mathcal{G})$ in general case depends from the genetic code $\mathcal{G}$
but not only from the group $G(\mathcal{G})$. For example,  $\mathcal{G}_{\varnothing}^l$ and
$\mathcal{G}_0^l$ define the isomorphic groups but
the corresponding braid groups $B_n(\mathcal{G}_{\varnothing}^l)$ and
$B_n(\mathcal{G}_{0}^l)$ are not isomorphic for $n > 1$.

Let's try to clarify how the corresponding braid groups
$B_n(\mathcal{G}^l_{\varnothing}),$ $B_n(\mathcal{G}^l_{0})$ and $B_n(\mathcal{G}^l_{1})$
are connected
for $n \geq 2$.

By the symbol $\mathcal{G}_{\varnothing }$ we will further denote  the genetic code
$\mathcal{G}^{\varnothing }_{\varnothing }$, or
genetic code $\mathcal{G}^{l}_{\varnothing }$. By the symbol $\mathcal{G}_{0}$
we will further denote  the genetic code $\mathcal{G}^{\varnothing }_{0}$,
or genetic code $\mathcal{G}^{l}_{0}$.
By the symbol $\mathcal{G}_{1}$ we will denote the genetic code $\mathcal{G}_{1}^l$, containing
a nontrivial relation. To indicate the explicit form of the relation
$r$ of the genetic code $\mathcal{G}^{l}_{1}$ we will write $\mathcal{G}^{l}_{1}(r)$.

Note that in the group
 $B_n(\mathcal{G}_{\varnothing})$ does not hold the relation (7).
Hence, there exists homomorphisms $\mu_1 : B_n(\mathcal{G}_{\varnothing}) \longrightarrow
B_n(\mathcal{G}_{0})$, $\mu_2 : B_n(\mathcal{G}_{\varnothing}) \longrightarrow
B_n(\mathcal{G}_{1})$. And the kernel of the first homomorphism
${\rm ker}(\mu_1)$ is the normal closure of the element
$$
\sigma_1 \sigma_2 \ldots \sigma_{n-2} \sigma_{n-1}^2 \sigma_{n-2} \ldots \sigma_2 \sigma_1
$$
in the group $B_n(\mathcal{G}_{\varnothing})$, while the kernel of the second homomorphism
is the normal closure of the element
$$
r^{-1} \sigma_1 \sigma_2 \ldots \sigma_{n-2} \sigma_{n-1}^2 \sigma_{n-2} \ldots \sigma_2
\sigma_1.
$$

For the group $B_n(\mathcal{G})$, $n > 1$, there is a homomorphism $\nu $ onto the symmetric
group $S_n$ which is defined on the generators by equations:
$$
\nu (\sigma_i) = (i, i+1),~~~i=1,2, \ldots, n-1;~~~~
\nu (x_j) = 1,~~~j=1,2, \ldots, l.
$$
The kernel of this homomorphism ${\rm ker}(\nu)$ we will call {\it the pure braid group} on
$n$ strands and denote by $P_n(\mathcal{G})$.

\medskip

To define the structure of $B_n(\mathcal{G})$ we define its subgroup
 $D_n = D_n(\mathcal{G})$ which includes such elements that under the action of homomorphism
$\nu$ map to permutations which fixing the symbol 1.

Let us introduce the following definitions:
$$
x_{1i} = x_i,~~~x_{k+1,i}=\sigma_k^{-1} x_{ki} \sigma_k,~~i=1,2,\ldots,
l,~~~1 \leq k \leq n-1.
$$
In this definitions the following lemma is true, and it is an analog of the corresponding lemma
from the paper of Zariski.

\medskip

{\bf Lemma 1.} {\it In the group $D_n$ the following formulas of conjugation are true}
$$
\begin{array}{lll}
1) & \sigma_{k}^{\varepsilon }x_{1i}\sigma_{k}^{-\varepsilon }=x_{1i},
~~~~~~~~~~~~~~\varepsilon = \pm 1,~~~ & k=2,3,\ldots,n-1;~~i=1,2,\ldots,l\\
& \\
2) & x_{2i}a_{12}x_{2i}^{-1}=x_{1i}^{-1}a_{12}x_{1i}, & \\
& \\
3) & x_{2i}^{-1}a_{12}x_{2i}=a_{12}^{-1}x_{1i}a_{12}x_{1i}^{-1}a_{12}, & \\
& \\
4) & x_{2i}^{-\varepsilon }a_{1j}x_{2i}^{\varepsilon }=a_{1j}, & j > 2, \\
& \\
5) & x_{2i}^{-1}x_{1i}x_{2i}=a_{12}^{-1}x_{1i}a_{12}, &  \\
& \\
6) & x_{2i}x_{1i}x_{2i}^{-1}=x_{1i}^{-1}a_{12}x_{1i}a_{12}^{-1}x_{1i}, &  \\
& \\
7) & x_{2i}^{-1}x_{1j}x_{2i}=a_{12}^{-1}x_{1j}, & i< j, \\
& \\
8) & x_{2i}x_{1j}x_{2i}^{-1}=(x_{1i}^{-1}a_{12}x_{1i})x_{1j}, & i < j,\\
& \\
9) & x_{2j}x_{1i}x_{2j}^{-1}=(x_{1j}^{-1}a_{12}x_{1j})a_{12}^{-1}x_{1i}
(x_{1j}^{-1}a_{12}^{-1}x_{1j}), & i < j,\\
& \\
10) & x_{2j}^{-1}x_{1i}x_{2j} = (a_{12}^{-1}x_{1j}a_{12}x_{1j}^{-1})x_{1i}a_{12}, &
i < j.\\
\end{array}
$$

\medskip

{\bf Proof.} The formula  1) follows from the defining relation (4).

From the relation  (5) follows the relation
$$
x_{1i}(\sigma_1^{-1}x_{1i}\sigma_1^{-1}) =
(\sigma_1^{-1}x_{1i}\sigma_1^{-1})x_{1i},~~~i=1,2,\ldots,l.
\eqno{(10)}
$$
From this
$$
(\sigma_1 x_{1i}^{-1} \sigma_1) x_{1i} (\sigma_1^{-1}x_{1i}\sigma_1^{-1})=x_{1i},~~~i=1,2,\ldots,l.
$$
Conjugating both hand sides of this equation by element $a_{12}=\sigma_1^2$ we get the equation
5). Then conjugating the equation 5) by element $x_{2i}^{-1}$ we get
$$
x_{1i} = x_{2i}a_{12}^{-1}x_{1i}a_{12}x_{2i}^{-1}.
\eqno{(11)}
$$
On the other hand from (10) follows that
$$
\sigma_1^{-1}x_{1i}\sigma_1^{-1}x_{1i} =
a_{12}^{-1}x_{1i}\sigma_1^{-1}x_{1i}\sigma_1,
$$
i.~e.
$$
(\sigma_1^{-1}x_{1i}\sigma_1)a_{12}^{-1}x_{1i} =
a_{12}^{-1}x_{1i}(\sigma_1^{-1}x_{1i}\sigma_1).
$$
Since
 $\sigma_1^{-1}x_{1i}\sigma_1 = x_{2i}$ we see that the elements $x_{2i}$ and
$a_{12}^{-1}x_{1i}$ are permutable.
Hence, we can write (11) in the form
$$
x_{1i} = a_{12}^{-1}x_{1i}x_{2i}a_{12}x_{2i}^{-1},
$$
i. e.
$$
x_{1i}^{-1}a_{12}x_{1i} = x_{2i}a_{12}x_{2i}^{-1}, \eqno{(12)}
$$
but this is the equation  2).

If we conjugate both hand sides of relation (12) by element $\sigma_1$ then we get
$$
(\sigma_1^{-1}x_{1i}^{-1}\sigma_1) a_{12}
(\sigma_1^{-1}x_{1i}\sigma_1) = a_{12}^{-1}x_{1i}a_{12}x_{1i}^{-1}a_{12}
$$
and the equation 3) is ascertained.

To prove 4) we remember that $a_{1j}=\sigma_1^{-1} \ldots \sigma_{j-2}^{-1}
\sigma_{j-1}^2\sigma_{j-2} \ldots \sigma_1$ and so
$$
x_{2i}^{\varepsilon } a_{1j} x_{2i}^{-\varepsilon }=
(\sigma_1^{-1}x_{1i}^{\varepsilon }\sigma_1) a_{1j}
(\sigma_1^{-1}x_{1i}^{-\varepsilon }\sigma_1) =\sigma_1^{-1} x_{1i}^{\varepsilon}a_{2j}
x_{1i}^{-\varepsilon } \sigma_1 = \sigma_{1}^{-1}a_{2j}\sigma_{1}
=a_{1j},
$$
i. e. the equation 4) is true.

Conjugate both hand sides of the equation 5) by element $a_{12}^{-1}x_{1i}$, we get
$$
x_{1i}^{-1}a_{12}x_{2i}^{-1}x_{1i}x_{2i}a_{12}^{-1}x_{1i} = x_{1i}.
$$
Using that element $x_{2i}$ is permutable with  $a_{12}^{-1}x_{1i}$ we rewrite
the last equation in the form
$$
x_{2i}^{-1}x_{1i}^{-1}a_{12}x_{1i}a_{12}^{-1}x_{1i}x_{2i} =
x_{1i}.
$$
Conjugate both hand sides of this relation by element $x_{2i}^{-1}$ we get the relation 6).

Observe  that the defining relation (6) of $B_n(\mathcal{G})$ we can write in the form
$$
x_{1j}x_{2i}^{-1}x_{1j}^{-1}x_{2i} = a_{12},~~~j > i.
$$
From this equation follows 7).

To prove 8) we conjugate both hand sides of equation 7) by element
$x_{2i}^{-1}$, we get
$$
x_{1j} = x_{2i}a_{12}^{-1}x_{1j}x_{2i}^{-1},~~~j
> i.
$$
Multiplying both hand sides of this equation on the left by $x_{2i}a_{12}x_{2i}^{-1}$ and
using
2) we get the required equation.

Conjugating both hand sides of the equation
$$
x_{1j}x_{2i}^{-1}x_{1j}^{-1}x_{2i} = a_{12},~~~j > i,
$$
which is true in the group $B_n(\mathcal{G})$ (see relation (6)) by the generator $\sigma_1$,
we get
$$
x_{2j}a_{12}^{-1}x_{1i}^{-1}a_{12}x_{2j}^{-1}a_{12}^{-1}x_{1i}a_{12} = a_{12}.
$$
Multiplying both hand sides of this equation on the right by
$a_{12}^{-1}x_{1i}^{-1}a_{12}$ and turning to inverse elements
we get the equation
$$
(x_{2j}a_{12}^{-1}x_{2j}^{-1})(x_{2j}x_{1i}x_{2j}^{-1})(x_{2j}a_{12}x_{2j}^{-1})=
a_{12}^{-1}x_{1i}.
$$
From this
$$
x_{2j}x_{1i}x_{2j}^{-1} = (x_{2j}a_{12}x_{2j}^{-1})a_{12}^{-1}x_{1i}(x_{2j}a_{12}^{-1}x_{2j}^{-1}).
$$
Using the equation 2) we get 9).

To prove 10) we conjugate both hand sides of 9)
by $x_{2j}$:
$$
x_{1i} = (x_{2j}^{-1}x_{1j}^{-1}x_{2j})(x_{2j}^{-1}a_{12}x_{2j})(x_{2j}^{-1}x_{1j}x_{2j})
(x_{2j}^{-1}a_{12}^{-1}x_{2j})(x_{2j}^{-1}x_{1i}x_{2j})\times
$$
$$
\times (x_{2j}^{-1}x_{1j}^{-1}x_{2j})
(x_{2j}^{-1}a_{12}^{-1}x_{2j})(x_{2j}^{-1}x_{1j}x_{2j}).
$$
Using the equations 3) and 5) we get the equation
$$
x_{1i} = x_{1j}a_{12}^{-1}x_{1j}^{-1}a_{12}(x_{2j}^{-1}x_{1i}x_{2j})a_{12}^{-1}.
$$
From this
$$
x_{2j}^{-1}x_{1i}x_{2j} =
(a_{12}^{-1}x_{1j}a_{12}x_{1j}^{-1})x_{1i}a_{12},
$$
but it is the equation 10). The lemma is proven.

\medskip

{\bf Lemma 2.} 1) {\it As a set of right coset representatives of $D_n$ in  $B_n(\mathcal{G})$
we can take the set
$M_n=\{1, \sigma_1, \sigma_2\sigma_1, \ldots , \sigma_{n-1}\sigma_{n-2} \ldots \sigma_2\sigma_1\}$.}
2) {\it The group $D_n$
is generated by elements} $a_{12}$, $a_{13}, \ldots, a_{1n},$
$\sigma_2, \sigma_3, \ldots, \sigma_{n-1},$ $x_{1k},$ $x_{2k},$ $k=1,2,\ldots,l$.
3) {\it The subgroup $U_n=\langle a_{12}, a_{13}, \ldots, a_{1n}, x_{11}, x_{12},\ldots, x_{1l}
\rangle $
is normal in} $D_n$.

\medskip

{\bf Proof.} 1) We denote by
$m_0=1,$ $m_i=\sigma_{i}\sigma_{i-1} \ldots \sigma_2\sigma_1,$ $i=1,2,\ldots,n-1$.
Note that the image of the element $m_i$ in
$S_n$ maps the symbol  $1$ to $i+1$ (we view the multiplication of permutations from right to
left). So $M_n$ is a set of right coset representatives of $D_n$ in $B_n(\mathcal{G})$.
And what is more,  $M_n$
is a Schreier set of coset representatives.

To prove  2) we use the Reidemeister-Schreier method. Let $u\longmapsto \overline{u}$
be a function which takes the right coset representatives of $D_n$ in $B_n(\mathcal{G})$. Then
$$
D_n=\langle d_{x,m} = (\overline{xm})^{-1}xm ~\vert ~m \in M_n, x \in
\{\sigma_{1}, \sigma_{2}, \ldots, \sigma_{n-1}, x_{1}, x_{2},\ldots,
x_{l}\} \rangle.
$$
We will show that if $x=\sigma_i$ for some  $i=1,2,\ldots,n-1$ then
$$
d_{\sigma_i,m_j}=
\left\{
\begin{array}{ll}
\sigma_i & \mbox{if } i > j+1, \\
1 & \mbox{if } i=j+1,\\
a_{1,i+1} & \mbox{if } i = j, \\
\sigma_{i+1} & \mbox{if } i < j.
\end{array} \right.
$$
Indeed, for $i > j+1$ we have
$$
d_{\sigma_i,m_j} = (\overline{\sigma_i m_j})^{-1} \sigma_i m_j = m_j^{-1}\sigma_i m_j =
\sigma_i.
$$
For $i=j+1$
$$
d_{\sigma_i,m_{i-1}} = (\overline{\sigma_i m_{i-1}})^{-1} \sigma_i m_{i-1} =
m_i^{-1}\sigma_i m_{i-1} = m_i^{-1}m_i=1.
$$
For $i=j$

$$
d_{\sigma_i,m_i} = (\overline{\sigma_i\sigma_i m_{i-1}})^{-1} \sigma_i^2 m_{i-1} =
m_{i-1}^{-1}\sigma_i^2 m_{i-1} =
m_{i-2}^{-1}\sigma_{i}\sigma_{i-1}^2\sigma_{i}^{-1}m_{i-2}=
$$

$$
=\sigma_i m_{i-3}^{-1}(\sigma_{i-2}^{-1}\sigma_{i-1}^2\sigma_{i-2})m_{i-3}\sigma_{i}^{-1}=
\sigma_i m_{i-3}^{-1} \sigma_{i-1} \sigma_{i-2}^2\sigma_{i-1}^{-1}m_{i-3}
\sigma_{i}^{-1}=
$$

$$
\sigma_i\sigma_{i-1}m_{i-3}^{-1}
\sigma_{i-2}^{2}m_{i-3}\sigma_{i-1}^{-1}\sigma_i = \ldots
=\sigma_{i}\sigma_{i-1} \ldots \sigma_2\sigma_{1}^2 \sigma_{2}^{-1} \ldots
\sigma_{i-1}^{-1}\sigma_{i}^{-1} = a_{1i}.
$$
At the end, for $i < j$ the requirement equation is got from the following
line:

$$
d_{\sigma_i,m_j} = (\overline{\sigma_i m_{j}})^{-1} \sigma_i m_{j} =
m_{j}^{-1}\sigma_i m_{j} =
$$

$$
=\sigma_{1}^{-1}\sigma_{2}^{-1} \ldots \sigma_{i+1}^{-1}
(\sigma_{i+2}^{-1}\ldots \sigma_{j}^{-1}\sigma_{i}\sigma_{j}\ldots
\sigma_{i+2})
\sigma_{i+1}\ldots \sigma_{2} \sigma_{1}=
$$

$$
=m_{i+1}^{-1}\sigma_{i}m_{i+1}=
m_i^{-1}\sigma_{i+1}^{-1}(\sigma_{i}\sigma_{i+1}\sigma_{i})m_{i-1}=
$$

$$
=m_i^{-1}\sigma_{i+1}^{-1}(\sigma_{i}\sigma_{i+1}\sigma_{i})m_{i-1}=
m_i^{-1}\sigma_{i+1}^{-1}\sigma_{i+1}\sigma_i\sigma_{i+1}m_{i-1}=
$$

$$
=m_i^{-1} \sigma_{i}m_{i-1}
\sigma_{i+1} = m_i^{-1}m_i\sigma_{i+1} = \sigma_{i+1}.
$$\\

Further, if $x=x_i$, $i=1,2,\ldots,l,$ then
$$
d_{x_i,m_j}=
\left\{
\begin{array}{ll}
x_i & \mbox{if } j=0, \\
x_{2i} & \mbox{if } j=1,2,\ldots,n-1.
\end{array} \right.
$$
Indeed,
$$
d_{x_i,1} = (\overline{x_{i}})^{-1} x_{i} = x_i,~~~
d_{x_i,m_j} = \overline{m_{j}}^{-1} x_{i}m_j = \sigma_1^{-1}x_i\sigma_1 =x_{2i},
$$
$$
~~~~~~~~~~~~~~~~~~~~~~~~~~~~~~~~~~~~~~~~~~~~~~ i=1,2,\ldots,l,~~~ j=1,2,\ldots,n-1.
$$
Hence, the item 2) of lemma is proven.

The item 3) immediately follows from the item 2), conjugation rules in $B_n$, and lemma 1.

\medskip

Using this lemma we can represent every element $w \in B_n(\mathcal{G} )$ as the product
$$
w=w_1m_{n,i_1},~~~w_1 \in D_n,~~~m_{n,i_1}= m_{i_1} = \sigma_{i_1} \sigma_{i_1-1} \ldots
\sigma_1 \in M_n,
$$
and the element $w_1$ in one's part we can represent in the form
$$
w_1=w_2 u_n,
$$
where $w_2 \in \langle \sigma_{2}, \sigma_{3}, \ldots ,
\sigma_{n-1}, x_{21}, x_{22}, \ldots ,x_{2l} \rangle$, $u_n \in U_n.$
Denote by $\overline{B}_{n-1}$ a subgroup of $B_n(\mathcal{G})$ which is generated by elements
$\sigma_{2}, \sigma_{3}, \ldots ,
\sigma_{n-1},$ $x_{21}, x_{22}, \ldots ,x_{2l}$. It is true

\medskip

{\bf Lemma 3.} {\it In the group $\overline{B}_{n-1}=\langle \sigma_{2}, \sigma_{3}, \ldots ,
\sigma_{n-1}, x_{21}, x_{22}, \ldots ,x_{2l}\rangle $ the following relations are true

$$
\sigma_i\sigma_{i+1}\sigma_i=\sigma_{i+1}\sigma_i\sigma_{i+1},~~~ i=2,3,...,n-2, \eqno{(13)}
$$

$$
\sigma_i \sigma_j = \sigma_j \sigma_i,~~~ |i-j|\geq 2,
\eqno{(14)}
$$

$$
\sigma_i x_{2k} = x_{2k} \sigma_i,~~~ i=3,4,...,n-1,~~~k=1,2,...,l,
\eqno{(15)}
$$

$$
x_{2k} (\sigma_2^{-1} x_{2k} \sigma_2^{-1}) = (\sigma_{2}^{-1} x_{2k} \sigma_2^{-1}) x_{2k},
~~~ k=1,2,...,l, \eqno{(16)}
$$

$$
x_{2j} (\sigma_2^{-1} x_{2i}^{-1} \sigma_2) x_{2j}^{-1}  (\sigma_{2}^{-1} x_{2i} \sigma_2)
= \sigma_{2}^2,
~~~ 1\leq i < j \leq l, \eqno{(17)}
$$

$$
r_2 = \sigma_2 \sigma_3 \ldots \sigma_{n-2} \sigma_{n-1}^{2} \sigma_{n-2} \ldots
\sigma_3 \sigma_{2} a_{12}.
 \eqno{(18)}
$$\\
where the word $r_2$ is got from $r_1=r$ conjugation by element $\sigma_1$, i. e.
if we replace all elements $x_{1k}=x_k$, which be included in $r$ by elements} $x_{2k}$,
$k=1,2,\ldots , l$.

\medskip

{\bf Proof.} Correctness of the relations (13)--(14) is evident.
The relation (15) is got from the relation (4) for $i=3,4,\ldots , n-1$ conjugation by the
generator $\sigma_1$. To prove the relation (16) we consider the relation (5):
$$
x_{1k} (\sigma_1^{-1} x_{1k} \sigma_1^{-1}) = (\sigma_{1}^{-1} x_{1k} \sigma_1^{-1}) x_{1k},
~~~ k=1,2,...,l,
$$
from $B_n(\mathcal{G})$. Conjugate its both hand sides by element $\sigma_2$ and
using permutability  $\sigma_2$ and $x_{1k}$ we get
$$
x_{1k} (\sigma_2^{-1} \sigma_1^{-1} \sigma_{2} x_{1k} \sigma_2^{-1} \sigma_1^{-1}
 \sigma_2) = ( \sigma_{2}^{-1} \sigma_1^{-1} \sigma_2 x_{1k} \sigma_2^{-1} \sigma_{1}^{-1}
 \sigma_2)  x_{1k}.
$$
Using the relation
$\sigma_2^{-1} \sigma_1^{-1} \sigma_{2} = \sigma_1 \sigma_2^{-1} \sigma_1^{-1}$ we get
$$
x_{1k} (\sigma_1 \sigma_2^{-1} \sigma_{1}^{-1} x_{1k} \sigma_1 \sigma_2^{-1} \sigma_1^{-1})
= (\sigma_{1} \sigma_2^{-1} \sigma_1^{-1} x_{1k} \sigma_1 \sigma_{2}^{-1} \sigma_1^{-1}) x_{1k}.
$$
Conjugating  both hand sides of this equation by element $\sigma_1$ we have
$$
(\sigma_1^{-1} x_{1k} \sigma_1) \sigma_2^{-1} (\sigma_{1}^{-1} x_{1k} \sigma_1) \sigma_2^{-1}
= \sigma_{2}^{-1} (\sigma_1^{-1} x_{1k} \sigma_1) \sigma_{2}^{-1} (\sigma_1^{-1}
x_{1k} \sigma_1).
$$
Since
$\sigma_1^{-1} x_{1k} \sigma_1 = x_{2k}$ then from this we get the relation (16).

To prove the relation (17) we consider the relation
$$
x_{1j} (\sigma_1^{-1} x_{1i}^{-1} \sigma_1) x_{1j}^{-1} (\sigma_1^{-1} x_{1i} \sigma_1) =
\sigma_{1}^{2},~~~j > i,
$$
which is true in $B_n(\mathcal{G} )$ (see the relation (6)). Conjugating both hand sides
of this relation by element $\sigma_2$, we get
$$
\sigma_2^{-1} x_{1j} (\sigma_1^{-1} x_{1i}^{-1} \sigma_1)
x_{1j}^{-1} (\sigma_{1}^{-1} x_{1i} \sigma_{1}) \sigma_2 = \sigma_2^{-1} \sigma_{1}^{2}
\sigma_2.
$$
Take advantage of permutability $\sigma_2$ with every generator
$x_{1k}$ we rewrite the last equality in the following form
$$
x_{1j} (\sigma_2^{-1} \sigma_1^{-1} \sigma_2) x_{1i}^{-1}
(\sigma_2^{-1}\sigma_1 \sigma_2)
x_{1j}^{-1} (\sigma_{2}^{-1} \sigma_1^{-1} \sigma_2) x_{1i} (\sigma_{2}^{-1} \sigma_1
 \sigma_2) = \sigma_{2}^{-1} \sigma_1^2 \sigma_2.
$$
Using the equations
$$
\sigma_2^{-1} \sigma_{1}^{-1} \sigma_2 = \sigma_1 \sigma_2^{-1} \sigma_{1}^{-1},~~~
\sigma_2^{-1} \sigma_1 \sigma_2 = \sigma_{1} \sigma_2 \sigma_1^{-1},
$$
we get
$$
x_{1j} (\sigma_1 \sigma_2^{-1} \sigma_1^{-1}) x_{1i}^{-1}
(\sigma_1 \sigma_2 \sigma_1^{-1})
x_{1j}^{-1} (\sigma_{1} \sigma_2^{-1} \sigma_1^{-1}) x_{1i} (\sigma_{1}
\sigma_2 \sigma_1^{-1}) = \sigma_{1} \sigma_2^2 \sigma_1^{-1}.
$$
Conjugating  both hand sides of this equation by element $\sigma_1$, and using that
$x_{2k} = \sigma_{1}^{-1}  x_{1k} \sigma_{1}$, $k=1, 2, \ldots, l$, we get the relation (17).

The relation (18) is obtained from the relation (7) conjugation by element $\sigma_1$.
The lemma is proven.

\medskip

From this lemma follows that the group $\overline{B}_{n-1}$ is isomorphic to the braid group
$B_{n-1}(\mathcal{G})$ in the case then there is not the relation
(7). Hence, in this case the following assertion is true

\medskip

{\bf Lemma 4.} {\it For the genetic code $\mathcal{G} = \mathcal{G}_{\varnothing }$
(doesn't contain relations)
the group $D_n$ contains a subgroup which isomorphic to $B_{n-1}(\mathcal{G})$ and so $D_n$
is the semi-direct product: $D_n=U_n \leftthreetimes B_{n-1}(\mathcal{G})$.}

\medskip

In the case when $\mathcal{G} = \mathcal{G}_1$ (contains only one non-trivial relation)
it is easy to see that the factor-group
 $D_n/U_n$ is isomorphic to  $B_{n-1}(\mathcal{G})$, but in this case we can't
assert that $B_{n-1}(\mathcal{G})$ is a subgroup of the group $B_{n}(\mathcal{G})$.

Further we can consider the group $B_{n-1}(\mathcal{G})=\langle \sigma_2, \ldots, \sigma_{n-1},
x_{21}, \ldots, x_{2l}\rangle, $ its homomorphism
onto the symmetric group $S_n$ which is induced by the homomorphism
$\nu : B_{n}(\mathcal{G})\longrightarrow S_n$, choose the subgroup
$D_{n-1} \leq B_{n-1}(\mathcal{G})$ that consists from the elements, which under
the action of $\nu $ map to the permutations that fix the symbol 2.
As a set of coset representatives of $D_{n-1}$ in $B_{n-1}(\mathcal{G})$ we can take the set
$$
M_{n-1}=\{1, \sigma_2, \sigma_3\sigma_2, \ldots, \sigma_{n-1} \sigma_{n-2} \ldots \sigma_2 \}.
$$
As above, we can show that $D_{n-1}$ is generated by elements
$$
a_{23}, a_{24}, \ldots, a_{2n}, \sigma_{3}, \sigma_{4}, \ldots, \sigma_{n-1},
x_{2k}, x_{3k}, ~~~k=1, 2, \ldots, l,
$$
and the subgroup
$U_{n-1} = \langle a_{23}, a_{24}, \ldots, a_{2n}, x_{21}, x_{22}, \ldots, x_{2l}\rangle $
is normal in
$D_{n-1}$. In this case the factor-group  $D_{n-1}/U_{n-1}$
is isomorphic to $B_{n-2}(\mathcal{G})$.
Continuing this process, on the last step we get a group  $B_{1}(\mathcal{G})$
which is isomorphic to $G(\mathcal{G})$.
In that way, for any  $i=0, 1, \ldots, n-2$ we construct a sequence
$$
B_{n-i}(\mathcal{G}) \geq D_{n-i} \unrhd U_{n-i},
$$
where
$$
B_{n-i}(\mathcal{G}) = \langle \sigma_{i+1}, \sigma_{i+2}, \ldots, \sigma_{n-1},
x_{i+1,1}, \ldots, x_{i+1,l} \rangle,
$$
$$
D_{n-i}
 = \langle a_{i+1,i+2}, a_{i+1,i+3}, \ldots, a_{i+1,n},
\sigma_{i+2}, \sigma_{i+3}, \ldots, \sigma_{n-1},
   x_{i+1,1}, \ldots,
$$
$$
~~~~~~~~~~~~~~~~~~~~~~~~~~~~~~~~~~~~~~~~~x_{i+1,l}, x_{i+2,1}, \ldots, x_{i+2,l} \rangle,
$$
$$
U_{n-i} = \langle a_{i+1,i+2}, a_{i+1,i+3}, \ldots, a_{i+1,n}, x_{i+1,1}, x_{i+1,2}, \ldots,
x_{i+1,l} \rangle.
$$
In this case the factor group $D_{n-i}/U_{n-i}$ is isomorphic to the group
$B_{n-i-1}(\mathcal{G})$.

Recall that the pure braid group $P_n(\mathcal{G})$ on $n$ strands was defined as
the kernel of the homomorphism $\nu : B_n(\mathcal{G})\longrightarrow S_n$.
>From the assertions which was proved above follows

\medskip

{\bf Lemma 5.} 1) {\it The subgroup $P_n(\mathcal{G})$ is a normal subgroup  of index  $n!$ in
the group} $B_n(\mathcal{G})$. 2) {\it As a set of right coset representatives of
$P_n(\mathcal{G})$ in  $B_n(\mathcal{G})$ we can take the set} $M = M_2 M_3 ... M_n.$
3) {\it The group $P_n(\mathcal{G})$ is generated by elements}
$$
a_{ij},~~~1\leq i < j \leq n;~~~x_{ik},~~~1\leq i \leq n,~~1\leq k \leq
l.
$$

\medskip

In the group $B_n(\mathcal{G})$ together with the relation (7) are contained relations
which getting from its by conjugation. More precisely, the following assertion is true

\medskip

{\bf Lemma 6.} {\it In the group $B_n(\mathcal{G})$ the following equations are true
$$
\begin{array}{l}
r_1 = a_{12} a_{13} \ldots a_{1n}, \\
r_2 = a_{23} a_{24} \ldots a_{2n} a_{12}, \\
r_3 = a_{34} a_{35} \ldots a_{3n} a_{13} a_{23}, \\
................................................\\
r_{n-1} = a_{n-1,n} a_{1,n-1} a_{2,n-1} \ldots a_{n-2,n-1}, \\
r_n = a_{1n} a_{2n} \ldots a_{n-1,n}, \\
\end{array}
$$
where the word $r_i$ is got from the word $r$ by replacement the generator
 $x_j$ by $x_{ij}$, $i=1,2,\ldots,n$, $j=1,2,\ldots,l$. Every relation
$$
r_i = a_{i,i+1} a_{i,i+2} \ldots a_{i,n}
$$
is true in the group  $U_{n-i+1}$, $i=1,2,\ldots,n-1$.}

\medskip

{\bf Proof.} Since  $x_k = x_{1k}$, $k=1,2,\ldots,l$ and
$$
\sigma_1 \sigma_2 \ldots \sigma_{n-2} \sigma_{n-1}^2 \sigma_{n-2} \ldots \sigma_2 \sigma_1 =
a_{12} a_{13} \ldots a_{1n}
$$
then the first relation coincides with the relation
(7) of $B_n(\mathcal{G})$. The relation
$$
r_i = a_{i,i+1} a_{i,i+2} \ldots a_{in} a_{1i} a_{2i} \ldots a_{i-1,i}
$$
is got from the first relation  using the conjugation  by element $\sigma_1 \sigma_2 \ldots
\sigma_{i-1}$.
Then we use the following conjugations rules:
$$
(\sigma_1 \sigma_2 \ldots \sigma_{i-1})^{-1} x_{1k} (\sigma_1 \sigma_2 \ldots
\sigma_{i-1}) = x_{ik},~~~k = 1, 2, \ldots, l,
$$
and the equations
$$
(a_{12} a_{13} \ldots a_{1n})^{\sigma_1 \sigma_2 \ldots \sigma_{i-1}} =
(\sigma_i \sigma_{i+1} \ldots \sigma_{n-1}^2 \ldots \sigma_{i+1} \sigma_i)
(\sigma_{i-1} \sigma_{i-2} \ldots \sigma_2 \sigma_{1}^2 \sigma_2 \ldots \sigma_{i-2}
\sigma_{i-1}) =
$$
$$
= (a_{i,i+1} a_{i,i+2} \ldots a_{i,n})(a_{1i} a_{2i} \ldots
a_{i-1,i}).
$$

Note that in the group
$$
U_{n} = \langle a_{12}, a_{13}, \ldots, a_{1n}, x_{11}, x_{12},\ldots,
x_{1l} \rangle
$$
the relation $r_1 = a_{12} a_{13} \ldots a_{1n}$ is true. Using the induction by $i$, and
equations of above,
we see that in the group
$$
U_{n-i+1} = \langle a_{i,i+1}, a_{i,i+2}, \ldots, a_{i,n}, x_{i1}, x_{i,2},\ldots, x_{il} \rangle
$$
only the following equation
$$
r_i = a_{i,i+1} a_{i,i+2} \ldots a_{i,n},~~~i=1,2,\ldots,n-1
$$
is true.
The lemma is proven.

\vskip 20pt

\begin{center}

{\bf \S~4. Some properties of  $B_n(\mathcal{G})$}

\end{center}

\vskip 5pt


We will show that $B_n(\mathcal{G}_0)$ for $n \geq 2$ has a torsion. It is true

\medskip

{\bf Theorem 1.} {\it The group $B_n(\mathcal{G}_0)$ has a torsion for all $n \geq 2$.
The subgroup $U_3(\mathcal{G}_0)$ has elements of order 2. The subgroups $U_n(\mathcal{G}_0)$ for
$n \neq 3$
are free.}

\medskip

{\bf Proof.}
Indeed, for $n=2$ in the group $B_n(\mathcal{G}_0)$ the following relation  $a_{12} = 1$ is true,
i.~e.
$\sigma_1^2=1$. If $\mathcal{G}_0 = \mathcal{G}_0^{\varnothing}$ then $B_2(\mathcal{G}_0) \simeq
\mathbb{Z}_2$.
If $\mathcal{G}_0 = \mathcal{G}_0^l$ then the group  $B_2(\mathcal{G}_0)$ is generated
by elements $x_1, x_2, \ldots, x_l, \sigma_1$. Let $x_{1k}=x_k$, $x_{2k} =
\sigma_1^{-1} x_{1k} \sigma_1 = \sigma_1 x_{1k} \sigma_1$, $k=1,2,\ldots,l$. Then relations
(5)--(7) have the following form:

$$
x_{1j} x_{2j} = x_{2j} x_{1j},~~~j = 1,2,\ldots,l,
$$

$$
x_{1j} x_{2i}^{-1} = x_{2i}^{-1} x_{1j},~~~j > i,
$$

$$
\sigma_1^2 = 1.
$$
From these relations follows that the pure braid group $P_2(\mathcal{G}_0^l) \simeq F_l \times
F_l$ is the direct product of two copies of the free group $F_l$. The group $U_2$
is generated by elements $x_{1k}$, $k=1,2,\ldots,l$, and so is free, and
$U_1 = \langle x_{21}, x_{22},\ldots,x_{2l} \rangle $ is free too.

By the lemma 6 in  $B_n(\mathcal{G}_0)$ for $n=3$ the following relations
$$
a_{12}a_{13} = 1,~~~a_{23}a_{12} = 1,~~~a_{13}a_{23} = 1
$$
are true.
If we express from the first relation $a_{12}$ and substitute it in the second then we get
$a_{23}a_{13}^{-1} = 1$, i. e. $a_{23} = a_{13}$. Then from the third relation we have
 $a_{23}^2 = 1$.
Hence, $a_{13} = a_{23},$ $a_{12} = a_{13}^{-1} = a_{23}^{-1} = a_{23}$.
Then
$$
U_3 = \langle a_{13}, x_{11}, x_{12},\ldots,x_{1l} \rangle \simeq
\mathbb{Z}_2 * F_l,
$$
$$
U_2 = \langle x_{21}, x_{22},\ldots,x_{2l} \rangle \simeq F_l,
$$
$$
U_1 = \langle x_{31}, x_{32},\ldots,x_{3l} \rangle \simeq F_l \simeq
G(\mathcal{G}_0).
$$

For $n=4$ in the group $B_n(\mathcal{G}_0)$ the following relations hold
$$
a_{12}a_{13}a_{14} = 1,~~~a_{23}a_{24}a_{12} = 1,~~~a_{34}a_{13}a_{23} = 1,
~~~a_{14}a_{24}a_{34} = 1.
$$
Let us express from the first relation $a_{14}$, from the second $a_{24}$, from the third
 $a_{34}$ and substitute these expressions in the last relation,
we get
$$
a_{14}a_{24}a_{34} = (a_{12}^{-1} a_{23}^{-1} a_{13}^{-1})^2 = 1,
$$
i. e.
$$
(a_{13}a_{23}a_{12})^2 = (a_{12} a_{13} a_{23})^2 = 1.
$$
Since $a_{12}a_{13}a_{23} = (\sigma_1 \sigma_2)^3$, and the image of the element $\sigma_1 \sigma_2$
in the group
$S_4$ different from 1 then the group
$B_4(\mathcal{G}_0)$ contains some elements of finite order.

Consider the subgroup
$$
U_4 =  \langle a_{12}, a_{13}, a_{14}, x_{11},
x_{12},\ldots,x_{1l} \rangle.
$$
In this group the relation $a_{12}a_{13}a_{14} = 1$ is true. If we exclude from it $a_{12}$,
then we see that
$$
U_4 = \langle a_{13}, a_{14}, x_{11},
x_{12},\ldots,x_{1l} \rangle
$$
is a free group. Further,
$$
U_3 = \langle a_{23}, a_{24},  x_{21},
x_{22},\ldots,x_{2l}~ \vert \vert ~a_{23}a_{24}=1 \rangle = \langle a_{24}, x_{21},
x_{22},\ldots,x_{2l} \rangle,
$$
but since $(a_{12}a_{13}a_{23})^2 = 1$ then in the group  $U_3$  the following
equation
$a_{23}^2 = a_{24}^{-2} = 1$ is true.
Hence, $U_3 \simeq \mathbb{Z}_2 * F_l$.

Consider
$$
U_2 = \langle a_{34},  x_{31},
x_{32},\ldots,x_{3l} \vert \vert a_{34}=1 \rangle =  \langle x_{31},
x_{32},\ldots,x_{3l} \rangle.
$$
Hence, $U_2$ as well as the group $U_1 = G(\mathcal{G}_0)$ is free.

Further, for arbitrary $n \geq 4$ in view of the lemma 6
the following equations
$$
\left.
\begin{array}{ll}
a_{12} a_{13} \ldots a_{1n} = 1, & \\
& \\
a_{23} a_{24} \ldots a_{2n} a_{12} = 1, &    \\
&\\
a_{34} a_{35} \ldots a_{3n} a_{13} a_{23} = 1, &    \\
& \\
....................................... & \\
& \\
a_{n-1,n} a_{1,n-1} a_{2,n-1}\ldots a_{n-2,n-1} = 1, &    \\
& \\
a_{1n} a_{2n} \ldots a_{n-1,n} = 1, &    \\
\end{array}
\right.
$$
are true in the group  $B_n(\mathcal{G}_0)$.
Express from the first equation $a_{1n}$, from the second $a_{2n}$, from the third
 $a_{3n}$, and so on, finally, from $n-1$-th we express $a_{n-1,n}$. If we substitute these
expressions in the last equation and use the defining relations of the group $P_n$ then we get
$$
\left[ a_{12}^{-1} (a_{23}^{-1} a_{13}^{-1}) (a_{34}^{-1} a_{24}^{-1}
a_{14}^{-1}) \ldots (a_{n-2,n-1}^{-1} a_{n-3,n-1}^{-1} \ldots a_{1,n-1}^{-1})\right]^2 =1.
$$
It is not hard to check that an element which is inverse to the element from the square bracket,
is equal to
$$
a_{12} (a_{13} a_{23}) (a_{14} a_{24} a_{34}) \ldots
(a_{1,n-1} a_{2,n-1} \ldots a_{n-2,n-1}),
$$
which is equal to
$$
(\sigma_1 \sigma_2 \ldots \sigma_{n-2} )^{n-1}
$$
(see \cite[c. 28]{Bir1}). Hence, $(\sigma_1 \sigma_2 \ldots \sigma_{n-2} )^{2n-2} = 1$.
It is easy to see that the image of the element
 $\sigma_1 \sigma_2 \ldots \sigma_{n-2} $ in the group
 $S_n$ different from 1, and so the group $B_n(\mathcal{G}_0)$ has elements of finite order.

Now consider the subgroups $U_{n-i+1},$ which have the genetic code
$$
U_{n-i+1} = \langle a_{i,i+1}, a_{i,i+2}, \ldots, a_{i,n} x_{i1},
x_{i2},\ldots,x_{il}~ \vert \vert ~a_{i,i+1} a_{i,i+2} \ldots a_{i,n} =
1 \rangle.
$$
If we exclude from the relation the generator $a_{i,i+1}$ then we see that
$$
U_{n-i+1} = \langle a_{i,i+2}, a_{i,i+3}, \ldots, a_{i,n}, x_{i1},
x_{i2},\ldots,x_{il} \rangle,~~~i = 1, 2, \ldots, n-1.
$$
For all $i\neq n-2$ this group is free, and for $i=n-2$ we have
$$
U_3 = \langle a_{n-2,n-1}, a_{n-2,n}, x_{n-2,1},
x_{n-2,2}, \ldots, x_{n-2,l}~ \vert \vert ~a_{n-2,n-1} a_{n-2,n} =
1,~~~a_{n-2}^2 = 1 \rangle \simeq \mathbb{Z}_2 * F_l.
$$
Thus the theorem is proven.

\medskip

This theorem generalizes a result from \cite{GilB} which asserts that the braid group of sphere
$B_n(S^2)$ has a torsion.

If we use the same argumentation as in  proof of the theorem  1 then we easy get

\medskip

{\bf Corollary.} {\it The subgroup $U_n(\mathcal{G}_1)$ of the braid group $B_n(\mathcal{G}_1)$,
$n \geq 2$,
is free.}

\medskip

Similar to the classical braid groups it is holds

\medskip

{\bf Theorem 2.} {\it In the group $B_n(\mathcal{G})$ the word problem is solable.}

\medskip

{\bf Proof.} Let a word $w$ be written in the generators of
$B_n(\mathcal{G})$. As was proved in the previous section there exists the short exact sequence
$$
1 \longrightarrow P_n(\mathcal{G}) \longrightarrow B_n(\mathcal{G}) \longrightarrow S_n
\longrightarrow 1.
$$
If the image of $w$ in $S_n$ is non-trivial then $w\neq 1$. In the opposite case,
$w \in P_n(\mathcal{G})$. Let us rewrite this word in the generators of the group
$P_n(\mathcal{G})$. Then we
represent $w$ as a product $w = w_1 w'$, where $w_1 \in U_n$ and is a word in generators
$a_{1j}$, $j=2,3, \ldots, n$, $x_{1k}$, $k=1,2, \ldots, l$, and the word $w'$ doesn't contain
these generators.
Now we consider the short exact sequence
$$
1 \longrightarrow U_n \longrightarrow P_n(\mathcal{G}) \longrightarrow
P_{n-1}(\mathcal{G})
\longrightarrow 1.
$$
The word $w_1 \in U_n$, and $U_n$ either a free group or is a free product with amalgamation,
and hence we can effectively check: is the word  $w_1$ represents trivial element. If $w_1\neq 1$
then $w\neq 1$. If $w_1 =1$ then we find an image of element $w'$ in the group
$P_{n-1}(\mathcal{G})$, and represent it as a product $w_2 w''$, where $w_2$ is a word
in the generators of $U_{n-1}$, and $w''$ doesn't contain these generators.
If  $w_2 \neq 1$ then  $w\neq 1$. If $w_2 = 1$ then we find the image of
 $w''$ in the group $P_{n-2}(\mathcal{G})$. If we continue this procedure then we either
prove that $w\neq 1$ or on a $n$-th step we construct an element $w_n \in U_1$.
Since  $U_1 \simeq G(\mathcal{G})$ and the word problem is effectively solvable in $G(\mathcal{G})$
then we can check: is $w_n = 1$ or not. If $w_n = 1$ then $w = 1$. In the opposite case
 $w\neq 1$. The theorem is proven.

\medskip

Recall, that the {\it rang } ${\rm rk}(G)$ of a finitely generated group $G$ \cite[p.
29]{OTG}
is called the minimal number of generators of $G$. For the trivial group
$G$ we let ${\rm rk}(G)=0$. It is holds

\medskip

{\bf Proposition 1.} {\it For the braid group $B_n(\mathcal{G})$ the following evaluations hold}
$$
\left.
\begin{array}{ll}
{\rm rk}(B_1(\mathcal{G})) = {\rm rk}(G(\mathcal{G})), &    \\
& \\
{\rm rk}(B_2(\mathcal{G})) \leq {\rm rk}(G(\mathcal{G})) + 1, &     \\
& \\
{\rm rk}(B_n(\mathcal{G})) \leq {\rm rk}(G(\mathcal{G})) + 2,~~~ \mbox{for }~~~n \geq 2.&
\end{array}
\right.
$$

\medskip

{\bf Proof.} The first equality follows from the definition since in this case
$B_1(\mathcal{G}) = G(\mathcal{G})$. For $n=2$ the group $B_2(\mathcal{G})$
is generated by element $\sigma_1$ and generators of $G(\mathcal{G})$, and hence the second
evaluation holds. The third evaluation follows from the fact that a subgroup of
 $B_n(\mathcal{G})$ which is generated by elements
$\sigma_1, \sigma_2, \ldots, \sigma_{n-1}$ can be generated by two elements  \cite[Chapter 6]{KM}.
The proposition is proven.

\medskip

In \cite{Bar-M92} was studied a question about the width of verbal subgroups
of the braid group $B_n(E^2)$ of Euclidean plane   $E^2$, and in \cite{Bar-T95} similar question
was studied for other Artin groups. Recall some facts from these papers.

The {\it width} ${\rm wid}(G, V)$ of verbal subgroup $V(G)$ which is defined in a group $G$
by a set of words $V$, with respect to this $V$ is called a minimal
$m \in \mathbb{N} \bigcup \{ +\infty \}$
such that every element from the subgroup $V(G)$ is a product of
 $\leq m$ values of words from $V \bigcup V^{-1}$. We will consider only finite sets of words
$V$ since for any verbal subgroup  $V(G)$ we can construct some set of words $W$ (which is infinite
in general case) such that $W(G) = V(G)$ and ${\rm wid}(G, W) =1$.

We will say that a group $G$ {\it is rich by verbal subgroups} if every set of words
$V$ that define a proper verbal subgroup of a free group $F_{2}$ rang two, define
a proper verbal subgroup of
$G$ (by the proper subgroup we mean a subgroup that is not trivial and is not all group).
In particular, if  $G$ contains a free non-abelian group then  $G$ is rich by verbal subgroups.

We will say that a group $G$ {\it doesn't have proper verbal subgroups
of finite width} if for every finite set $V$ that defines a proper verbal subgroup
 $V(G)$ of $G$
the width ${\rm wid}(G, V)$ is infinite.

Following statements hold

\medskip

{\bf Lemma 7}  (\cite[Lemma 1]{Bar-T95}). {\it Let $G=G_1 \times G_2 \times \ldots \times G_k$
be a direct product of finite number subgroups
 $G_i$. Then the width ${\rm wid}(G, V)$ is equal to the maximal width of the
factors ${\rm wid}(G_i, V)$ for all} $i=1, 2, \ldots, k$.

\medskip

{\bf Lemma 8}  (\cite[Lemma 3]{Bar-T95}). {\it If there is an epimorphism of  $G$ onto some group
which doesn't have proper verbal subgroups of finite width
and which is rich by verbal subgroups then  $G$ doesn't have proper verbal subgroups of
finite width.}

\medskip

Recall, that the braid group $B_n(\mathcal{G}^{1}_{\varnothing})$ is Artin group of type
 $B$. From some results of the papers
 \cite{Bar-M92, Bar-T95} follows that the braid group
$B_n(E^2)=B_n(\mathcal{G}^{\varnothing}_{\varnothing})$ of Euclidean plane, the braid groups
$B_n(\mathcal{G}^{1}_{\varnothing})$, and  corresponding pure braid groups $P_n(E^2)$ and
$P_n(\mathcal{G}^{1}_{\varnothing})$ for $n\geq 3$ doesn't have proper verbal subgroups of finite
width.

We suppose that holds

\medskip

{\bf Hypothesis.}  For every genetic code  $\mathcal{G}$ there is a natural number
$N = N(\mathcal{G})$ such that for all $n \geq N$ the braid group $B_n(\mathcal{G})$
doesn't have proper verbal subgroups of finite
width.

\medskip

For the pure  braid group holds

\medskip

{\bf Proposition 2.} {\it Let $\mathcal{G} = \mathcal{G}_{\varnothing}^l, \mathcal{G}_0^l$ for
$l\geq 2$ or
$\mathcal{G} = \mathcal{G}_{1}^l$ for $l\geq 3$. Then for every $n\geq 1$
the pure braid group $P_n(\mathcal{G})$
doesn't have proper verbal subgroups of finite
width.}

\medskip

{\bf Proof.}
Let $V$ be a finite proper set of words.
We can see that for the pure braid group $P_n(\mathcal{G})$ there is an epimorphism on $n$-th
direct power
$H=G(\mathcal{G}) \times G(\mathcal{G})\times \ldots \times G(\mathcal{G})$.
In  view of the lemma 7 ${\rm wid}(H, V) = {\rm wid}(G(\mathcal{G}), V)$. The group
$G(\mathcal{G})$ is either a free non-abelian or a group with one relation and  $\geq 3$
generators. In the both cases  $G(\mathcal{G})$  is rich by verbal subgroups.
The fact that $G(\mathcal{G})$ doesn't have proper verbal subgroups of finite
width follows from \cite[Theorem 2]{Bar-A97}. Hence,
in view of the lemma 8 the group  $P_n(\mathcal{G})$ also doesn't have proper verbal subgroups
of finite width.
The proposition is proven.

In connection with the defined construction we formulate the following questions.

\medskip

{\bf Question 1.} (I. A. Dynnikov). Let $\mathcal{G}$ be a genetic code with one
defining relation, $K_{\mathcal{G}}$ be an one-cell complex that corresponding
to the genetic code  $K_{\mathcal{G}}$, $B_n(K_{\mathcal{G}})$ be the braid group on $n$ strands,
that corresponding to the complex $K_{\mathcal{G}}$ (it is defined similar to the
braid group of manifold). In what cases  the braid group
$B_n(\mathcal{G})$ of the genetic code $\mathcal{G}$ is isomorphic to the braid group
$B_n(K_{\mathcal{G}})$ of the complex $K_{\mathcal{G}}$?

\medskip

{\bf Question 2.} Let  $M$ be a two dimensions compact non-orientable manifold of genus
 $g$, $\mathcal{G} = \langle x_1, x_2, \ldots, x_g~ \vert \vert ~x^2_1 x^2_2
\ldots x^2_g = 1)$
be a genetic code that defines the fundamental group $\pi_1(M)$.
How do  the groups $B_n(M)$ and $B_n(\mathcal{G})$ connect for $n \geq 2$?

\medskip

{\bf Question 3.} Which groups $B_n(\mathcal{G})$ are decomposeble by non-trivial manner in a
free product with amalgamation?

\end{document}